# An Amendment of Fast Subspace Tracking Methods

Zhu Cheng[1,2], Zhan Wang, Haitao Liu
1. School of Electron. Sci. & Eng.,
Nat. Univ. of Defense Technol.,
Changsha, China
Chengzhu88@gmail.com

Majid Ahmadi
2. Department of ECE,
University of Windsor, Windsor,
Canada ON
Ahmadi@uwindsor.ca

*Abstract*—**Tuning stepsize between convergence rate and steady state error level or stability is a problem in some subspace tracking schemes. Methods in DPM and OJA class may show sparks in their steady state error sometimes, even with a rather small stepsize. By a study on the schemes' updating formula, it is found that the update only happens in a specific plane but not all the subspace basis. Through an analysis on relationship between the vectors in that plane, an amendment as needed is made on the algorithm routine to fix the problem by constricting the stepsize at every update step. The simulation confirms elimination of the sparks.**

*Array signal processing; subspace tracking; stepsize; convergence.*

## I. INTRODUCTION

Tracking a subspace is Estimating a projection matrix onto a space or a basis for the space, from a random vector sequence observed by a sensor array. It is a powerful tool in some signal processing fields such as: telecommunication, radar, sonar and navigation, serves as a measure of adaptive filter, DOA estimation, or interference mitigation. Subspace tracking methodology is classified in into two categories: the first is estimating the space where the signal is generated from, the second one is to find orthogonal complement of that space. The former is known as a principal subspace (PS, PSA) tracker or signal subspace tracker, the later is often referred to as minor subspace (MS, MSA) tracker or a noise subspace tracker. For our earlier works on MUSIC, we are used to the term signal subspace track or noise subspace track.

N.L. Owsley developed the first algorithm for subspace tracking in [1]. Assuming the problem's dimension is N, the rank of the subspace we are interested in is L .Usually L << N. Complex of this solution proportion to $N^2L$,or $O(N^2L)$ namely. Many schemes with less compute complex were developed after then. An excellent survey paper [20] outlined almost all of achievements on this topic before 1990, which cost $O(N^2L)$ or $O(NL^2)$ operations. Algorithms with $O(NL)$ operations were developed after it. The new class is called as Fast Subspace Tracking method. Surveys on fast subspace schemes are presented in [3 pp30–43] or [19, pp221-270].

Let $x(k)$ is an N-dim observer vector from an N-element sensor array, as (1),

$$x(k) = \sum_{i=1}^{L} a_i s_i(k) + n(k), \qquad (1)$$

Where $a_i$ is N-dim vector with unit length, independent to each other, representing the manifold of one of the arriving signals, and $s_i(k)$ is a random variable independent to each other, representing the arriving signals from different source, $n(k)$ is a N-dim i.i.d random vector represent the sensor noises. Assuming $V = \{v_1, v_2, v_3 ..... v_L\}$ is a orthonormal basis of span( $a_i$ ). The signal subspace tracking is seeking a basis $W(k) = \{w_1(k), w_2(k), ..... w_L(k)\}$ which spans a subspace same as the space spanned by V; the noise subspace tracking is searching a $W(k)$ which spans a subspace as orthogonal compensate of span(V). From now, we may use the name of the basis as the name of the space, such as space W or space W(k) to span(W(k)), if there is no confusion.

One of criterion on subspace problem is the distance between the spaces. Majority of those solutions[2-11] use the projection error power as (2) for signal subspace tracking or (3) for noise subspace tracking.

$$e_p = \|P_V - P_W\|_F^2 \qquad (2)$$

$$e_p = \|P_{V\perp} - P_W\|_F^2 \qquad (3)$$

Where $P_V = VV^H$, $P_{V\perp} = I_N - P_V$, $P_W = WW^H$;

$\|\bullet\|_F$ means Frobenius Norm.

$I_N$ means N-dim identity matrix.

An additional criterion for orthonormality of W(k) as (4):

$$\eta = \|W^H(k)W(k) - I_L\|_F^2. \qquad (4)$$

DPM[21] class algorithm is started from optimization the coordinate length of input vector's image projected onto the subspace as (5) , while Oja class scheme optimize the length of the projection image as (6).

$$E(\|W^H x(k)\|^2) \qquad (5)$$

$$E(\|x(k) - P_W x(k)\|^2) \qquad (6)$$

The final routines of DPM or OJA type approaches are very similar. A typical routine of DPM class scheme is similar to (7).

This work was partly supported by National Natural Science Foundation of China under Grant(11173068))

$$q(k) = W^H(k-1)x(k)$$
$$y(k) = W(k-1)q(k)$$
$$T(k) = W(k-1) \pm \beta x(k)q^H(k) \quad (7)$$
$$W(k) = \text{orthnorm}\{T(k)\}$$

orthnorm(•) means orthonormalization operation.
Plus sign stands for signal subspace tracking,
Minus sign is for noise subspace tracking.

The typical routine for an Oja scheme only replaces $x(k)$ with $p(k) = x(k) - y(k)$ in the temporary general basis $T(k)$ update in (7). The variety of DPM might include FDPM, FRANS, HFRANS[6], MFPDM[7] and a version of SOOJA[8]. The branches of Oja include OOJA[9], OOJAH[9], FOOJA[10], original version of SOOJA[11].

Some unreasonable random sparks were observed when we apply these schemes, especially when the noise subspace tracking was applied. It happens to DPM class schemes under noise subspace tracking and all variety of Oja methods.

In this paper, we present the geometric relationship among x(k), y(k), the old last estimated space W(k-1) and the next estimated space W(k) by update equations analysis. An Amendment as needed on the schemes is made to fix the sparks problem by the applying a limiter on stepsize at every update step. The presented simulation confirmed the amendment.

## II. ANALYSIS OF THE UPDATE EQUATIONS

Considering the subspace tracking problem, the new arriving x(k) is projected onto the last estimated space W(k-1) to get a y(k) in space W(k-1). Vector t(k) is the y(k)'s projection image in x(k)'s orthogonal complement as (8). Certainly $t \perp x$.

$$t(k) = y(k) - x(k)x^H(k)y(k)/(x^H(k)x(k)) \quad (8)$$

The relation between those vectors and the last estimated space W(k-1) shows in Fig.1 with an omission of time index k. We omit the time index k in some of following equations if there is no confusion.

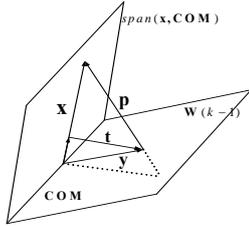

Figure 1. Relationship of vector and space under discussion

Assuming there is a vector set with L elements including y direction served as a orthonormal basis set for space W(k-1), the basis is noted by $(y/\|y\|, COM)$, where $COM$ is an L-1 dim subspace of space $W(k-1)$, and $COM \perp y(k)$. W(k-1) is another orthonormal basis for the same space. Both $(y/\|y\|, COM)$ and W(k-1) are orthonormal base set for the same subspace. Therefore an orthonormal matrix Q exists and meet (9,10);

$$(\mathbf{y}/\|\mathbf{y}\|, \mathbf{COM}) = W(k-1)\mathbf{Q}, \quad (9)$$
$$W(k-1) = (\mathbf{y}/\|\mathbf{y}\|, \mathbf{COM})\mathbf{Q}^H. \quad (10)$$

By left multiplication $W^H(k-1)$ to (9), it is easy to find the first row of Q is $q_0 = W^H(k-1)\frac{y}{\|y\|} = \frac{q}{\|y\|} = \frac{q}{\|q\|}$, therefore Q might be noted as $(q_0, \mathbf{Q}_{COM})$, with $q_0 \perp \mathbf{Q}_{COM}$.

$$x^H COM = x^H P_{W(k-1)} COM = x^H P_{W(k-1)}^H COM$$
$$= (P_{W(k-1)}x)^H COM = y^H COM = O \quad (11)$$

For (11) x is orthogonal to COM, therefore any linear combination of x and y is orthogonal to COM too, so is p.

Before the orthnorm(•) in (7), $T(k)$ is a general basis of the newly estimated subspace. Right multiply the update equation of $T(k)$ With $Q = (q_0, \mathbf{Q}_{COM})$. $T(k)Q$ is another base vector set for the newly estimated subspace as $Q = (q_0, \mathbf{Q}_{COM})$ is a unitary matrix.

For Oja schemes:
$$\mathbf{T}(k)Q = (\mathbf{w}(k-1) \pm \beta\mathbf{p}(k)q^H(k))(q_0, \mathbf{Q}_{COM})$$
$$= \mathbf{w}(k-1)(q_0, \mathbf{Q}_{COM}) \pm \beta\mathbf{p}(k)q^H(k)(q_0, \mathbf{Q}_{COM})$$
$$= (\mathbf{w}(k-1)q_0, \mathbf{w}(k-1)\mathbf{Q}_{COM})$$
$$\quad \pm \beta(\mathbf{p}(k)q^H(k)q_0, \mathbf{p}(k)q^H(k)\mathbf{Q}_{COM})$$
$$= (y/\|y\|, COM) \pm \beta(\mathbf{p}(k)\|q(k)\|, 0, 0, ..0)$$
$$= (y/\|y\| \pm \beta\mathbf{p}(k)\|q(k)\|, COM) \quad (12)$$
$$\triangleq (h_1, COM)$$

Or for DPM class:
$$\mathbf{T}(k)Q = (y/\|y\| \pm \beta x(k)\|q(k)\|, COM) \triangleq (h_2, COM). \quad (13)$$

After the multiplication, all newly estimated subspaces in Oja or DPM class schemes have a basis in term of COM and a linear combinations of x and y. By comparison on newly estimated subspaces basis (12) or (13) with another predefined basis $(y/\|y\|, COM)$ for the last estimated subspaces at k-1 time, it can be found that in the update routine COM keeps still while the basis vector updates only in the plane span(x,y). Fig.2 shows the update happening in span(x,y) as a cross-section of Fig1.

From Fig.2, the difference between signal or noise subspace tracking method is only the move direction of new basis vector h from the old base vector y to or from x. In span(x,y), t is furthest vector from x in angle view, and the closest vector to x is itself.

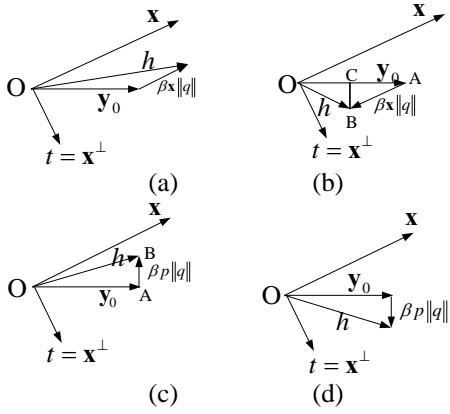

Figure 2 the relation of new candidate replace vector for y in span(x,y)

a) DPM class for signal subspace b) DPM class for Noise subspace

c) OJA class for signal subspace d) OJA class for Noise subspace

Therefore it is enough for the new vector h to move only between x and y for signal subspace tracking, outside that range cannot provide a better result than inside; for noise subspace tracking the boundary is between y and t. If the h is out of these ranges, sparks on the steady state projection error power might happen.

The relation between h and y or x is controlled by the step size parameter Beta. It is necessary to find a boundary for it.

### III. BOUNDARY OF BETA

From Fig.2a, for DPM type algorithm, in signal subspace tracking, any positive value of Beta always keep the h between x and y, therefore any positive Beta value will not cause the problem of stability. Therefore we can not find the sparks in DPM signal subspace tracking scheme.

But for DPM noise subspace tracking or all OJA classes, a fixed Beta value may move new base vector h beyond the predefined boundary x or t in section 2, when the input is big enough, it might causes deviate or spark. To avoid this situation, we try to find the boundary of Beta for them.

Statement: To avoid sparks, in noise subspace tracking, condition for DPM class algorithm is $\beta \le 1/\|x\|^2$, condition for Oja class schemes is $\beta \le 1/\|p\|^2$; in signal subspace tracking, condition for OJA class algorithm is $\beta \le 1/\|y\|^2$. Or generally $\beta \le 1/\|x\|^2$ is sufficient to avoid sparks under all those situations.

Proof: For DPM type noise subspace algorithm, with (13) and Fig2.b.

AO is vector $y/\|y\|$, therefore $|AO|=1$

AB is vector $\beta x\|q\|$, $\frac{y}{\|y\|}$ so $|AB| = \beta\|x\|\|q\| = \beta\|x\|\|y\|$

If the angle between x and y is $\theta$, then

$|BC| = \beta\|x\|\|y\|\sin\theta$, $|AC| = \beta\|x\|\|y\|\cos\theta$

therefore $|OC| = 1 - \beta\|x\|\|y\|\cos\theta$

If the angle between y and new vector h is $\gamma$, then

$$\tan\gamma = \frac{\beta\|x\|\|y\|\sin\theta}{1-\beta\|x\|\|y\|\cos\theta}.$$

For we set new vector not to over turn than the vector t, and the angle between y and t is $90^0 - \theta$, so $\gamma \le 90^0 - \theta$. For all the angles here range between 0 and 90deg, so:

$\tan\gamma \le \tan(90^0 - \theta)$. It means $\tan\gamma \le \cos\theta/\sin\theta$.

$$\frac{\beta\|x\|\|y\|\sin\theta}{1-\beta\|x\|\|y\|\cos\theta} \le \frac{\cos\theta}{\sin\theta}.$$

$$\beta \le \cos\theta/(\|x\|\|y\|) = 1/\|x\|^2$$

Therefore for FDPM type noise subspace tracking if $\beta \le 1/\|x\|^2$, then no overshoot will happen, and the sparks will cease.

Proof for OJA type scheme is similar to that of DPM algorithm. We omit it to save space.

S.Attallah and his colleagues had found a similar boundary from a different aspect for some of their OJA and FRANS schemes, to search a more aggressive update stepsize and improve convergent rate [15-18]. Unfortunately, there was no stable and fast subspace tracking scheme as FDPM or FOOJA for them to be used as a prototype at that time. From our view, aggressive stepsize is not our object, the stability of the scheme is. We are arming an amendment which can avoid the over tune phenomenon or eliminate the sparks.

### IV. AMENDMENT AND SIMULATION

From the Statement in section 3, if we set Beta as minima of all possible $1/\|x\|^2$, the sparks will cease by itself, the convergence rate will be rather slow. If we apply beta as $1/\|x(k)\|^2$ at any step, it will be too aggressive and make the steady state error rather high when the input is small.

We decide to use a reasonable Beta as the original algorithms, but at any time if $\beta \le 1/\|x\|^2$, we set $\beta = 1/\|x\|^2$ for that step.

We present parts of our simulations in order to verify the amendment. The simulation setup is similar to [2, 3], but differs in value of parameter. We consider a signal plus noise model with N =8, the random signal x(n) lies in an L=4 dimensional linear subspace, for convenience, assume manifolds of arriving signals $a_i$ are orthogonal to each other to avoid the interaction between the arrival angle and convergence rate, arriving signals from different source $s_i(k)$ are random variables independent to each other with the powers equal to [10, 1, 0.1, 0.1]. The noise $n(k)$ is N-dim iid white and Gaussian random vector, with $10^{-3}$ variance. Beta=0.08 for all simulation. Duration of simulation is 6000 steps. For each single run, there is a break of the basis at 3000 step to introduce more projection error and destroy the orthonormality, by adding of random matrix on them, every element of it is a iid variable with 0.1 variance, to check the ability of orthonormal and projection error power convergence. Only the results of FOOJA are displayed here, and results on FPDM or both version of SOOJA are similar to them.

Projection error power (3) or (4) is applied for comparison in Fig.3. The y axis is in db scale. From Fig.3, average of projection error power for the original version(blue without mark) is higher than the amended one(with+), and maxima projection error power out of 100 runs for the original versions(red without mark) is much higher and noise than that of the amended ones(black, with o).

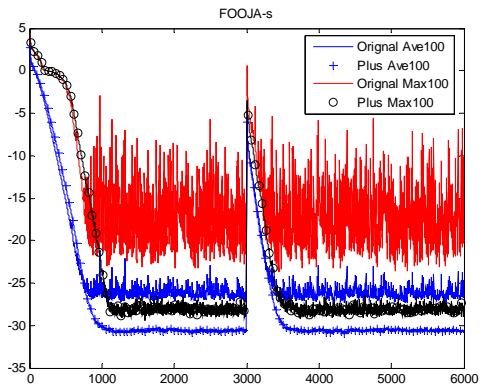

(a) FOOJA Signal Subspace

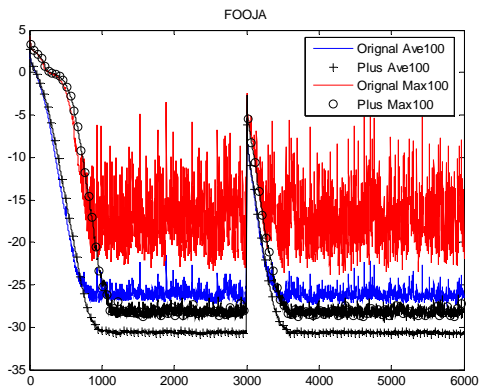

(b) FOOJA Noise Subspace

Figure 3 Compare the Amended FOOJA and original ones

## V. CONCLUSION

It is recommended to check and $\beta \leq 1/\|x\|^2$ apply a limiter on each iteration step for all mentioned schemes on noise subspace tracking, and OJA type schemes on signal subspace tracking.

REFERENCES

[1] N.L. Owsley, Adaptive Data Orthogonalization, ICASSP' 1978, pp. 09-112,197

[2] X.G.Doukopoulos,G.V.Moustakides, Fast and Stable Subspace Tracking, IEEE Trans on Signal Processing, VOL. 56, NO. 4, APRIL 2008 pp1452

[3] X.G.Doukopoulos, Power techniques for blind channel estimation in wireless communication systems, Ph.D.dissertation,IRISA-INRIA,Univ. Rennes, Rennes ,France, 2004.

[4] S.Attallah and K.Abed-Maraim, "Low-cost adaptive algorithm for noise subspace estimation," IEE Electronics Letters, vol. 38, no. 12, pp. 609-611, June 2002.

[5] E.Oja,"Principal components, minor components, and linear neural networks,"Neural Networks, vol. 5, pp. 927–935, Nov./Dec. 1992.

[6] S.Attallah,"The generalized Rayleigh's quotient adaptive noise subspace algorithm: A Householder transformation-based implementation," IEEE Trans. Circuits Syst. II, Exp. Briefs, vol. 42, no. 9, pp. 3-7, Jan. 2006

[7] Qian Linjie, Cheng Zhu, et al; Modification Algorithms for a class of subspace tracking methods; Signal Processing(Chinese); V26 Num5, May 2010, pp741

[8] Rong Wang, et al; A Novel Orthonormalization Matrix Based Fast and Stable DPM Algorithm for Principal and Minor Subspace Tracking; IEEE Trans on Signal Processing VOL. 60, NO.1, JAN. 2012.pp466

[9] K. Abed-Meraim, S. Attallah, A. Chkeif, and Y. Hua, "Orthogonal oja algorithm," IEEE Signal Processing Letters, pp. 116–119, 2000

[10] S. Bartelmaos and K. Abed-Meraim, "Principal and minor subspace tracking: Algorithms & stability analysis," in Proc. IEEE ICASSP'06, Toulouse, France, May 2006, vol. III, pp. 560-563

[11] Rong Wang; et al; Stable and Orthonormal OJA Algorithm With Low Complexity, IEEE Signal Processing Letters, Volume: 18 , Issue: 4 Year: 2011 , Page(s): 211 - 214

[12] Xiangyu, K.; Changhua, H.; Han, C.; A dual purpose principal and minor subspace gradient flow; IEEE Trans on Signal Processing, VOL. 60, NO. 1, JANUARY 2012,pp197

[13] Carl D. Meyer ;Matrix Analysis and Applied Linear Algebra;SIAM February 15, 2001

[14] Rong Wang, Minli Yao, Zhu Cheng, Hongxing Zou: Interference cancellation in GPS receiver using noise subspace tracking algorithm. Signal Processing 91(2): 338-343 (2011)

[15] Jonathan H.Manton,Iven Y.Mareeks,Samir Attallah;An analysis of the fast subspace tracking algorithm NOja; ICASSP,2002 Volume: 2 , pp1101-1104

[16] Samir Attallah and Karim Abed-Meraim,Fast Algorithms for Subspace Tracking, IEEE Signal Processing Letters, VOL. 8, NO. 7, JULY 2001 pp203

[17] Samir Attallah,REVISITING ADAPTIVE SIGNAL SUBSPACE ESTIMATION BASED ON RAYLEIGH'S QUOTIENT, ICASSP '01 Page(s): 4009 - 4012 vol.6

[18] Attallah, S.; Abed-Meraim, K.;Subspace estimation and tracking using enhanced versions of Oja's algorithm ;(SPAWC '01). Page(s): 218 – 220

[19] J.P. Delmas, "Ch.4 Subspace tracking for signal processing," Adaptive Signal Processing : Next Generation Solutions, Editors: T. Adali and S. Haykin, Wiley Interscience, April 2010.

[20] Comon, P.; Golub, G.H.; Tracking a few extreme singular values and vectors in signal processing;Proceedings of the IEEE,Volume: 78 , Issue: 8; Page(s): 1327 – 1343

[21] J. F. Yang and M. Kaveh, "Adaptive eigensubspace algorithms for directionor frequency estimation and tracking," IEEE Trans. Acoust. Speech Signal Processing, vol. 36, pp. 241–251, Feb. 1988.